\documentclass{amsart}

\usepackage{amsmath}
\usepackage{amssymb}
\usepackage{amsthm}
\usepackage{comment}
\usepackage{hyperref}
\usepackage{tikz}
\usetikzlibrary{cd}

\newcommand{\pol}{\mathrm{pol}}
\newcommand{\Quad}{\mathrm{Quad}}
\newcommand{\Mod}{\mathrm{Mod}}
\newcommand{\QQ}{\mathrm{Q}}

\newtheorem{theorem}{Theorem}[section]
\newtheorem{lemma}[theorem]{Lemma}
\newtheorem{definition}[theorem]{Definition}
\newtheorem*{udefinition}{Definition}
\newtheorem{proposition}[theorem]{Proposition}
\newtheorem{corollary}[theorem]{Corollary}
\newtheorem*{ucorollary}{Corollary}
\newtheorem*{utheoremA}{Theorem A}
\newtheorem*{utheoremB}{Theorem B}
\newtheorem*{uquestion}{Question}
\theoremstyle{definition}

\newtheorem{remark}[theorem]{Remark}

\title{A note on quadratic forms}

\author{Fabian Hebestreit}
\address{Department of Mathematics, University of Aberdeen, Scotland}
\email{fabian.hebestreit@abdn.ac.uk, hebestreit@math.uni-bielefeld.de}

\author{Achim Krause}
\address{Mathematisches Institut, Universität M\"unster, Germany}
\email{krauseac@uni-muenster.de}

\author{Maxime Ramzi}
\address{Institut for Matematiske Fag, Københavns Universitet, Denmark}
\email{maxime.ramzi@math.ku.dk}

\swapnumbers
\begin{document}

\begin{abstract}
For a field extension $L/K$ we consider maps that are quadratic over $L$ but whose polarisation is only bilinear over $K$. Our main result is that all such are automatically quadratic forms over $L$ in the usual sense if and only if $L/K$ is formally unramified. In particular, this shows that over finite and number fields, one of the axioms in the standard definition of quadratic forms is superfluous.
\end{abstract}

\maketitle
\setcounter{tocdepth}{1}
\tableofcontents
\section{Introduction}

Let $S$ be a commutative ring and $M, N$ two $S$-modules. Recall that an \emph{$S$-quadratic map} from $M$ to $N$ is a map $q: M\to N$ which satisfies
\begin{enumerate}
  \item $q(\lambda x) = \lambda^2 q(x)$
\item $q(x+y+z) + q(x) + q(y) + q(z) = q(x+y) + q(x+z) + q(y+z)$, and
\item $q(\lambda x + \mu y) + \lambda\mu(q(x)+q(y)) = \lambda\mu q(x+y) + q(\lambda x) + q(\mu y)$
%  \item $\pol_q(x+y,z) = \pol_q(x,z) + \pol_q(y,z)$,
%  \item $\pol_q(\lambda x,\mu y) = \lambda\mu \pol_q(x,y)$
\end{enumerate}
for all $\lambda,\mu \in S$ and $x,y,z \in M$. This note arose from the following question: \\

How independent are these three axioms? \\

Recall immediately, that the second and third axiom can be reformulated in terms of the polarisation 
\[\pol_q \colon M \times M \longrightarrow N, \quad (x,y) \longmapsto q(x+y) - q(x) - q(y),\]
with the second axiom equivalent to  $\pol_q(x+y,z) = \pol_q(x,z) + \pol_q(y,z)$ and the third to $\pol_q(\lambda x,\mu y) = \lambda\mu \pol_q(x,y)$, i.e.\ in total they are equivalent to the $S$-bilinearity of $\pol_q$. Their independence then has two easy cases: Every linear map $q$ has $\pol_q = 0$, so satisfies the second and third axioms, but the first axiom holds if and only $\lambda^2q(x) = \lambda q(x)$, which is not usually true, e.g. it is false for the identity map $S \rightarrow S$ of an integral domain other than $\mathbb F_2$. Secondly, for $S = \mathbb F_2$ any map with $q(0) = 0$ satifies the first and third axioms, but rarely the second, e.g. take the map $\mathbb F_2^3 \rightarrow \mathbb F_2$ which takes every non-zero vector to $1$. 

The remaining case seems more subtle, however: Since additive maps are linear over $\mathbb Z$, the second axiom implies the third for $S = \mathbb Z$ and this implication is preserved under epimorphisms of rings: Given the second condition, the third  can be rephrased (for $\mu = 1$, but this clearly suffices by symmetry of $\pol_q$) as the two ring homomorphisms $S \longrightarrow \mathrm{End}_\mathbb Z(B_q)$
\[\lambda \longmapsto [b \mapsto \lambda \cdot b] \quad \text{and} \quad \lambda \longmapsto [b \mapsto b(\lambda -,-)]\]
agreeing, where $B_q$ is the $S$-submodule of biadditive maps $M \times M \rightarrow N$ generated by the $\pol_q(\lambda-,-)$ for $\lambda \in S$. We thus conclude that the second axiom implies the third also for all quotients and localisations of $\mathbb Z$, such as the $\mathbb F_l$ and $\mathbb Q$. 

In a different direction, already the first axiom clearly implies the other two if $M$ is free of rank $1$. We were, however, very surprised to find that the first two axioms still imply the third for maps $\mathbb F_{p^2}^2 \rightarrow \mathbb F_{p^2}$ and $\mathbb Q[i]^2 \rightarrow \mathbb Q[i]$. Our main result is a generalisation of this observation. To state it let us introduce some notation:

\begin{udefinition}
  For a homomorphism $\varphi \colon R \rightarrow S$ among commutative rings and two $S$-modules $M$ and $N$, let $\Quad_\varphi(M,N)$ (or $\Quad_{S/R}(M,N)$) denote the $S$-module of maps $q \colon M \rightarrow N$ with
\begin{enumerate}
\item $q(s m) = s^2 q(m)$ for all $s\in S$,
\item $\pol_q \colon \varphi^*M \times \varphi^*M \rightarrow \varphi^*N$ is $R$-bilinear.
\end{enumerate}
We will refer to such maps as $\varphi$- or \emph{$S/R$-quadratic}, and abbreviate $\Quad_{\mathrm{id}_S}(M,N)$ to $\Quad_S(M,N)$.
\end{udefinition}

Note that $\Quad_{S/\mathbb Z}(M,N)$ is precisely the module of maps $M \rightarrow N$ which satisfy the first two axioms of an $S$-quadratic form as above, but not necessarily the third, and if $L$ is any of the prime fields, then $\Quad_{L/\mathbb Z}(M,N) = \Quad_L(M,N)$ for all $L$-vector spaces $M,N$. But more:

\begin{utheoremA}
  \label{thm:thmB}
  If $L/K$ is a field extension,  then $\Quad_{L/K}(M,N) = \Quad_L(M,N)$ for all $L$-vector spaces $M$ and $N$ if and only if the Kähler differentials $\Omega_{L/K}$ vanish, that is if and only if
\begin{enumerate}
\item the extension $L/K$ is generated by $p$-th powers for $\mathrm{char}(K) = p>0$, or
\item the extension $L/K$ is algebraic for $\mathrm{char}(K) = 0$.
\end{enumerate}
\end{utheoremA}

Note also straight from the definition that if $\Quad_{S/R}(M,N) = \Quad_{S}(M,N)$, then $ \Quad_{T/R}(M,N) = \Quad_{T/S}(M,N)$ for two homomorphisms $\varphi \colon R \rightarrow S$ and $\psi \colon S \rightarrow T$ and $T$-modules $M$ and $N$. We thus for example find the (at least to us) surprising statement:

\begin{ucorollary}
The third axiom in the definition of an $L$-quadratic form is superfluous whenever $L$ is an algebraic field extension of a prime field, e.g. when $L$ is a finite or a number field.
\end{ucorollary}

The connection to Kähler differentials in Theorem A comes from the fact that any $R$-linear derivation $d\colon S \rightarrow M$ gives an element $q_d \in \Quad_{\varphi}(S^2,M)$ via
\[q_d(s,s') = d(s)s'-sd(s').\]
For example, the derivative of polynomials gives the $R[T]/R$-quadratic form
\[q \colon R[T]^2 \rightarrow R[T], \quad (F,G) \longmapsto F'G-FG'\]
for $S = R[T]$ with 
\[ \pol_q\left((T,0),(0,1)\right) = 1 \neq 0 = T\pol_q((1,0),(0,1)),\]
which therefore witnesses that the inclusion \[\Quad_{R[T]}(R[T]^2,R[T]) \subseteq \Quad_{R[T]/R}(R[T]^2,R[T])\] is strict. For an example over fields, simply replace the polynomial ring by its fraction field. For an algebraic example (necessarily inseparable) note that the above form also defines an element of 
$\Quad_{R[T]/R[T^p]}(R[T]^2,R[T])$ whenever $R$ has characteristic $p$. In the final section of the text, we shall show that this construction accounts for all exotic forms in these cases.

In fact, if on the contrary $\Quad_{S/R}(S^2,N) = \Quad_{S}(S^2,N)$ for some $S$-module $N$, we have $\Quad_{S/R}(M,N) = \Quad_{S}(M,N)$ for all $S$-modules $M$: It suffices to observe that we can restrict a $\varphi$-quadratic form $q$ along any $S$-linear map $f:S^2\to M$, and that by definition, if all the $q\circ f$ are $S$-quadratic, then so is $q$. It therefore suffices to describe the difference between $\Quad_\varphi(S^2,N)$ and $\Quad_S(S^2,N)$ for all $S$-modules  $N$ and an arbitrary $\varphi: R \rightarrow S$. We do so by considering the $S$-algebra $\QQ_\varphi = (S \otimes_R S) \otimes_{\Delta_\varphi} S$ augmented over $S$ by multiplication, where $\Delta_\varphi$ is the $R$-subalgebra of $S \otimes_R S$ spanned by the elements $s \otimes s$. As an $S$-module $\QQ_\varphi$ then canonically splits as $\mathrm W_\varphi \oplus S$ (with $W_\varphi$ the kernel of the canonical map $\QQ_\varphi \to S$), and essentially by unwinding definitions,
one finds
\[\Quad_\varphi(S^2,N) \cong \Quad_S(S^2,N) \oplus \mathrm{Hom}_S(\mathrm W_\varphi,N),\]
so that $\Quad_\varphi(M,N) = \Quad_S(M,N)$ for all $S$-modules $M$ and $N$ if and only if $\mathrm W_\varphi = 0$ or equivalently $\QQ_\varphi \cong S$ via the canonical maps (see Lemma \ref{lm:quadsplit} and Proposition \ref{prop:criterionintro} for a precise discussion).

Analysing the ring $\QQ_\varphi$ further, we then find the following general criteria at and away from the prime $2$:

\begin{utheoremB}
For a ring homomorphism $\varphi \colon R \rightarrow S$ we have $\Quad_\varphi(M,N) = \Quad_S(M,N)$ for all $S$-modules $M$ and $N$ if and only if
\begin{enumerate}
  \item the map $\mathrm{Fr} \cdot \varphi \colon S \otimes R \rightarrow S$, or equivalently the relative Frobenius $\mathrm{Fr}_\varphi\colon S\otimes_R \mathrm{Fr}^* (R/2)\to S$ of $S/R$, is an epimorphism in case $2=0$ in $S$, or
\item the kernel of the multiplication map $S \otimes_R S \rightarrow S$ is an idempotent ideal, or equivalently $\Omega_{S/R} = 0$, in case $2$ is a unit in $S$.
\end{enumerate}
\end{utheoremB}

In fact, the construction of $\varphi$-quadratic forms from derivations provides a map $\mathrm W_\varphi \rightarrow \Omega_{\varphi}$, restricted from the map
\[\QQ_\varphi \longrightarrow \Omega_\varphi, \quad s\otimes s' \otimes t \longmapsto t(d(s)s' - sd(s')),\] 
which is always surjective, and which we show to be an isomorphism whenever $2$ is a unit in $S$. We thus obtain
\[\Quad_{S/R}(S^2,N) \cong \Quad_S(S^2,N)\oplus \mathrm{Der}_R(S,N)\]
away from characteristic $2$. In characteristic $2$ the connection between quadratic forms and derivations remains somewhat more mysterious (to us). Here we identify $\mathrm W_\varphi$ as the kernel of the multiplication map
\[\mathrm {Fr}_\varphi^* S \otimes_{\mathrm{Fr}^* R \otimes_R S} \mathrm {Fr}_\varphi^* S \longrightarrow S.\]
leading to the second condition above, but the map $\mathrm W_{S/R} \longrightarrow \Omega_{S/R}$ is not generally injective (see below). Nevertheless, we did not find an example in which the target vanishes but the source does not and indeed among fields no such example can exist, which leads to Theorem A. However, we do not know about the general case:

\begin{uquestion}
Is there a ring map $\varphi \colon R \rightarrow S$ in positive characteristic, such that $\Omega_{S/R} = 0$ without $\mathrm{Fr} \cdot \varphi \colon S \otimes R \rightarrow S$, or equivalently the relative Frobenius $\mathrm{Fr}_\varphi$ of $S/R$, being an epimorphism? 
\end{uquestion}

Slightly more drastically, one may ask for such an example already in the case $R = \mathbb F_p$, i.e.\ for a formally unramified $\mathbb F_p$-algebra $S$ whose Frobenius is not an epimorphism. Let us note immediately that examples where the Frobenius is not surjective (i.e.\ where $S$ is not semi-perfect) do indeed exist, whereas for $S$ of finite type surjectivity of Frobenius follows, see for example the discussion in \cite{MathOver}.

The simplest example we found of an $S/R$-quadratic form that does not arise from a derivation as above is given by 
\[\mathbb F_2[X,Y]^2 \longrightarrow \mathbb F_2[X,Y], \quad (F,G) \longmapsto \frac{\partial^2 FG}{\partial X\cdot \partial Y}\]
for $S/R = \mathbb F_2[X,Y]/\mathbb F_2$ and similarly for the quotient field.

\subsubsection*{Acknowledgements}
We heartily thank Julius Frank for an afternoon of working through the example $\mathbb Q[i]^2 \rightarrow \mathbb Q[i]$, Lukas Brantner for a very useful exchange about infinite purely inseparable extensions, Andy Senger for a discussion about the relative Frobenius, and Thomas Nikolaus for hosting FH and MR during a visit to Münster, where this note took shape.

AK was supported by by the German Research Foundation (DFG) via the collaborative research centre ``Geometry:\ Deformations and Rigidity" and the cluster ``Mathematics Münster:\ Dynamics--Geometry--Structure" under grant nos.\ SFB 1442--427320536 and   
EXC 2044--390685587, respectively. MR was supported by the Danish National Research Foundation (DNRF) through the ``Copenhagen Center for Geometry and Topology" under grant no.\ DNRF151.

\section{Preliminaries}

We start by considering the functor $\Quad_\varphi(-,N)$ for fixed $N$.

\begin{proposition}\label{propsifted}
The functor $\Quad_\varphi(-,N) \colon \Mod_S^{\mathrm{op}} \rightarrow \Mod_S$ preserves $1$-cosifted limits. In particular, for a free resolution $F \rightarrow M[0]$ over $S$, we have
\[\Quad_\varphi(M,N) = \mathrm{ker}[\Quad_\varphi(F_0,N) \xrightarrow{(\mathrm{pr}_2 + d_1\mathrm{pr_1})^* - \mathrm{pr}_2^*}\Quad_\varphi(F_1 \oplus F_0,N)].\]
and for a free $S$-module $F$ we have
\[\Quad_\varphi(F,N) = \lim_{G \subseteq M} \Quad_\varphi(G,N)\]
where the limit runs over the (opposite of the) poset of finitely generated free summands $G$ of $M$.
\end{proposition}

Recall that a category $I$ is called $1$-sifted if it is nonempty and its diagonal functor is $1$-cofinal, i.e.\ the colimit of every diagram $I \times I \rightarrow \mathrm{Set}$ agrees with that of the $I$-shaped diagram $I \rightarrow \mathrm{Set}$ obtained by precomposition with the diagonal of $I$; $I$ is called $1$-cosifted if $I^\mathrm{op}$ is $1$-sifted.

\begin{proof}
By design $1$-siftedness of $I$ implies that the functor $\mathrm{Set} \rightarrow \mathrm{Set}, X \mapsto X \times X$ preserves $I$-shaped colimits, which in turn implies the same also for the functor $X \mapsto S \times X^2$ and the like. This allows one to provide the structure of an $S$-module to the colimit in $\mathrm{Set}$ of a $1$-sifted diagram of $S$-modules, which then serves as a colimit in $S$-modules. More succinctly, the forgetful functor $\Mod_S \rightarrow \mathrm{Set}$ preserves $1$-sifted colimits.  

We can now write $\Quad_\varphi(M,N)$ as three iterated equalisers, each of which will preserve $1$-cosifted limits. Namely, define $Q_1(M)$ as the equaliser of the two maps
\[\mathrm{Hom}_{\mathrm{Set}}(M,N) \longrightarrow \mathrm{Hom}_{\mathrm{Set}}(S\times M,N)\]
given by
\[q \longmapsto [(s,m) \mapsto s^2 q(m)] \quad \text{and} \quad q \longmapsto [(s,m) \mapsto q(sm)].\]
Similarly, define $Q_2(M)$ as the equaliser of the two maps
\[Q_1(M) \longrightarrow \mathrm{Hom}_{\mathrm{Set}}(M^3,N)\]
given by
\[q \longmapsto [(x,y,z) \mapsto \pol_q(x+y,z)] \quad \text{and} \quad q \longmapsto [(x,y,z) \mapsto \pol_q(x,z)+\pol_q(y,z)]\]
and $Q_3(M)$ as the equaliser of the two maps
\[Q_2(M) \longrightarrow \mathrm{Hom}_{\mathrm{Set}}(S\times M^2,N)\]
\[q \longmapsto [(s,x,y) \mapsto \pol_q(sx,y)] \quad \text{and} \quad q \longmapsto [(s,x,y) \mapsto s\pol_q(x,y)].\]
  Then we have $Q_3(M) = \Quad_\varphi(M,N)$ by construction, and by the discussion at the beginning of the proof, all terms occuring in the construction preserve $1$-cosifted limits (in $M \in \Mod_S^{\mathrm{op}}$). Limits commuting with limits then implies the first claim.

  To see the other statements, note that for a free resolution $F \rightarrow M[0]$ we find that $M$ is the coequaliser of the two maps $\mathrm{pr_2}, \mathrm{pr}_2 + d_1\mathrm{pr_1} \colon F_1 \oplus F_0 \rightarrow F_0$, which admit $(0,\mathrm{id}_{F_0})$ as a common split (this is the degree $\leq 1$ part of the Dold-Kan correspondence). Since reflexive coequalisers are $1$-sifted it follows that $\Quad_\varphi(M,N)$ is the equaliser of the two maps 
\[\mathrm{pr_2}^*, (\mathrm{pr}_2 + d_1\mathrm{pr_1})^* \colon \Quad_\varphi(F_0,N) \longrightarrow\ \Quad_\varphi(F_1 \oplus F_0,N)\]
which is the kernel from the statement. The final statement follows from filtered categories being $1$-sifted.
\end{proof}

To study $\Quad_\varphi(-,N)$ we may thus restrict it to finitely generated free modules. Observe next that
\[\Quad_\varphi(S,N) = \Quad_S(S,N)\]
since both sides consist precisely of those maps $q \colon S \rightarrow N$ with $q(s) = s^2q(1)$; in particular, both sides are canonically isomorphic to $N$ via evaluation at $1$. To get from $M=S$ to an arbitrary finitely generated free module, we have:

\begin{lemma}\label{lm:quadsplit}
Let $\Delta_\varphi \subseteq S \otimes_R S$ the sub-$R$-algebra spanned by the elements of the form $s\otimes s$. Then $\Quad_\varphi(M \oplus M',N)$ is isomorphic to 
\[\Quad_\varphi(M,N) \ \oplus\ \mathrm{Hom}_S(S \otimes_{\Delta_\varphi} (M \otimes_R M'),N)\ \oplus\ \Quad_\varphi(M',N),\]
where the outer two maps are induced by the restrictions to $M$ and $M'$, respectively, and the middle map is induced by restricting the polarisation along $M \otimes M' \xrightarrow{\mathrm{incl}_1 \otimes \mathrm{incl}_2} (M \oplus M') \otimes (M \oplus M')$. In particular, the decomposition in natural in all three variables.
\end{lemma}

\begin{remark}
In particular, $\Quad_\varphi(-,N)$ is a quadratic functor. Since it also preserves $1$-cosifted limits, it is classified by a form parameter as in \cite{Schlichting} or \cite[Section 4.2]{P1}. This is given by
\[\mathrm{Hom}_{\Delta_\varphi}(S \otimes_R S,N) \xrightarrow{\mathrm{ev}_{1 \otimes 1}} N \xrightarrow{2\mu^*} \mathrm{Hom}_{\Delta_\varphi}(S \otimes_R S,N),\]
where $\mathrm{C}_2$ acts on source and target by flipping the $S$ factors and the middle term $N$ is acted on by $S$ via the squaring map on $S$.
\end{remark}

\begin{proof}
By general non-sense the extra summand in $\Quad_\varphi(M \oplus M',N)$ is exactly the kernel of the outer two restriction maps, i.e.\ it consists of those quadratic maps $M \oplus M' \rightarrow N$ that vanish on $M \oplus 0$ and $0 \oplus M'$. But this means that $q(m,m') = \pol_q((m,0),(0,m'))$. This observation provides a bijection between such quadratic maps and $R$-bilinear maps $b \colon M \times M' \rightarrow N$, which satisfy $b(sm,sm') = s^2b(m,m')$ for all $s \in S$. But these correspond to precisely those maps $M \otimes_R M' \rightarrow N$ that are linear over the multiplication map $\Delta_\varphi \rightarrow S$, and this is equivalent to the claim.
\end{proof}

\begin{definition}
We shall denote the $S$-algebra $(S \otimes_R S) \otimes_{\Delta_\varphi} S$ by $\QQ_\varphi$.
\end{definition}

Here, the $S$-algebra structure is given by $s \mapsto 1 \otimes 1 \otimes s$. Note that $\QQ_\varphi$ is also augmented over $S$ via the natural map $\QQ_\varphi \rightarrow \QQ_{S/S} = S$ induced by the natural inclusion $\Quad_S(S^2,-) \subseteq \Quad_\varphi(S^2,-)$. In particular,
\[\QQ_\varphi \cong S \oplus \mathrm{W}_\varphi\]
as $S$-modules, with $\mathrm{W}_\varphi$ the kernel of the augmentation. Unwinding definitions this augmentation is simply given by the multiplication on $S$. 

In total, we thus obtain the general criterion mentioned in the introduction, and in particular the statement that the inclusion $\Quad_S(M,N) \subseteq \Quad_\varphi(M,N)$ is an equality whenever $\varphi$ is an epimorphism of rings:

\begin{proposition}\label{prop:criterionintro}
For an $S$-module $N$, the inclusion $\Quad_S(M,N) \subseteq \Quad_\varphi(M,N)$ is an equality for all $S$-modules $M$ if and only if it is one for $M=S^2$ if and only if the multiplication of $S$ induces an isomorphism
\[N = \mathrm{Hom}_S(S,N) \longrightarrow \mathrm{Hom}_S(\QQ_\varphi,N),\]
and generally
\[\Quad_\varphi(S^2,N) \cong \Quad_S(S^2,N) \oplus \mathrm{Hom}_S(\mathrm{W}_\varphi,N).\]
\end{proposition}

\begin{proof}
The first statement is immediate from the previous lemma, and the second one follows by observing that under the decompositions
\[\Quad_S(S^2,N) = \Quad_S(S,N) \ \oplus\ \mathrm{Hom}_S(\QQ_{S/S},N)\ \oplus\ \Quad_S(S,N)\]
and 
\[\Quad_\varphi(S^2,N) = \Quad_\varphi(S,N) \ \oplus\ \mathrm{Hom}_S(\QQ_\varphi,N)\ \oplus\ \Quad_\varphi(S,N)\]
the inclusion $\Quad_S(S^2,N) \subseteq \Quad_\varphi(S^2,N)$ is induced by by the identity on the outer summands and the augmentation on the middle summand.
\end{proof}

As the final preparation we note:

\begin{lemma}\label{lemma:flatfixed}
If $\varphi \colon R \rightarrow S$ makes $S$ a flat $R$-module, then $\Delta_\varphi = (S \otimes_R S)^{\mathrm C_2}$, where the action is by flipping the factors.
\end{lemma}

\begin{proof}
Clearly, $\Delta_\varphi \subseteq (S \otimes_R S)^{\mathrm C_2}$ holds for arbitrary $\varphi \colon R \rightarrow S$. For the converse consider the statement that for an $R$-module $M$ the $R$-module $(M \otimes_R M)^{\mathrm C_2}$ is generated by the diagonal terms $m \otimes m$. Then this is clearly true for $M=R$. But 
\[((M \oplus M')\otimes_R (M \oplus M'))^{\mathrm C_2} \cong (M \otimes_R M)^{\mathrm C_2} \oplus (M \otimes M') \oplus (M' \otimes_R M')^{\mathrm C_2}\]
with $(m \otimes m, a \otimes b, m' \otimes m')$ on the right mapping to 
  $m\otimes m + m'\otimes m' +(a+b)\otimes (a+b) - a\otimes a - b\otimes b$
  on the left, so that if the statement is true for $M$ and $M'$ then it is also true for $M\oplus M'$, and similarly if the statement is true for some filtered system of $R$-modules, then it also true for their colimit, since filtered colimits commute with fixed points. It follows that the statement is in particular true for any filtered colimit of free modules and these are precisely the flat ones by Lazard's theorem.
\end{proof}

\begin{remark}
This result is certainly not optimal: For example, the conclusion is also true for every epimorphism $R \rightarrow S$. By the classification of finitely generated modules it then follows that $\Delta_\varphi = (S \otimes_R S)^{\mathrm C_2}$ for every algebra $S$ over a principal ideal domain $R$. The statement also holds whenever $2$ is a unit in $S$, since then the norm map $(S \otimes_R S)_{\mathrm C_2} \rightarrow (S \otimes_R S)^{\mathrm C_2}$ is surjective and one easily checks that it always takes values in $\Delta_\varphi$. 

The conclusion does not, however, hold in complete generality. For a simple counterexample take $R = \mathbb Z[X,Y,Z]$, let $I$ be its augmentation ideal and put $S = R \oplus I$, the associated square zero extension, with $\varphi$ the evident inclusion. Then the element $XY \otimes Z \in I \otimes_R I$ is a fixed point for the flip-action, but we claim it does not lie in the submodule $\Delta_I$ of $I \otimes_R I$ generated by the diagonal elements, and the same then follows for its image in $S \otimes_R S$. To see this consider the image under the multiplication $I \otimes_R I \rightarrow R/2$. The image of $\Delta_I$ is the ideal $(X^2,Y^2,Z^2)$, which clearly does not contain $XYZ$.
\end{remark}

\section{Proof of the main results}

In the present section we shall deduce the main results of this note. We start with the case away from $2$:

\begin{proposition}\label{ident2inv}
Let $I$ be the kernel of the multiplication $S \otimes_R S \rightarrow S$. If $2$ is a unit in $S$, then the map 
\[\iota \colon S \longrightarrow (S \otimes_R S)/I^2, \quad s \longmapsto \frac{1 \otimes s + s \otimes 1} 2\]
is a ring homomorphism and the two maps
\[\QQ_\varphi \longrightarrow (S \otimes_R S)/I^2, \quad s \otimes s' \otimes t \mapsto \frac{st \otimes s' + s \otimes s't} 2\]
\[(S \otimes_R S)/I^2 \longrightarrow \QQ_\varphi, \quad s \otimes s' \longmapsto s \otimes s' \otimes 1\]
are inverse isomorphisms over and under $S$.
\end{proposition}

In other words, whenever $2$ is a unit in $S$ the map 
\[S^2 \longrightarrow (S \otimes_R S)/I^2, \quad (s,t) \longmapsto s \otimes t\]
is $S/R$-quadratic, where we give the target the module structure from above, and is in fact the universal such, that vanishes on $S \times 0$ and $0 \times S$.

Before we dive into the proof we record the following calculation, which we will reuse in the proof of Proposition \ref{ident2=0} below: In $(S \otimes_R S) \otimes_{\Delta_\varphi} S$ we have
\begin{align*}
s \otimes s' \otimes 1 + s' \otimes s \otimes 1 &= (s+s') \otimes (s+s') \otimes 1 - s \otimes s \otimes 1 - s'\otimes s' \otimes 1 \\
&= 1 \otimes 1 \otimes [(s+s')^2-s^2-s'^2] \label{eq}\tag{$\ast$}\\
&= 1 \otimes 1 \otimes 2ss' 
\end{align*} 
for every $\varphi \colon R \rightarrow S$ and $s,s' \in S$,

\begin{proof}
We first check that $\iota$ is multiplicative. To this end we observe
\[\frac{1 \otimes st + st \otimes 1} 2 - \frac{1 \otimes s + s \otimes 1}2\cdot \frac{1 \otimes t + t \otimes 1}2 = \frac{(1 \otimes s - s \otimes 1)(1 \otimes t - t \otimes 1)} 4 \]
which is in $I^2$. Next, we check that the two maps in question are well-defined. For the first one, we observe that
\[\frac{st \otimes s' + s \otimes s't} 2 = (s \otimes s')\iota(t)\]
and 
\[\iota(u)^2 - u \otimes u = \frac{(1 \otimes u - u \otimes 1)^2} 4 \in I^2\]
so that 
\[(s \otimes s')(u \otimes u)\iota(t) = (s \otimes s')\iota(u^2)\iota(t) = (s \otimes s')\iota(u^2t)\]
which together with bilinearity over $R$ implies bilinearity over $\Delta_\varphi$. For the second map note that for any $x = \sum_i s_i \otimes s_i' \in S \otimes_R S$, there is a (unique) decomposition $x = x_s + x_a$ with $x_s \in S \otimes_R S$ fixed by the flip of the two factors, and $x_a \in I$ acquiring a sign upon flipping, i.e. if $x= \sum_i s_i \otimes s'_i$ then 
\[x_s = \sum_i \frac{s_i \otimes s'_i + s'_i \otimes s_i} 2 \quad \text{and} \quad x_a = \sum_i \frac{s_i \otimes s'_i - s'_i \otimes s_i} 2.\]
If now $x \in I$, then by (\ref{eq}) above $x_s$ is mapped to 
\[\sum_i \frac{s_i \otimes s'_i \otimes 1 + s'_i \otimes s_i \otimes 1} 2 = 1 \otimes 1 \otimes \sum_i{s_is_i'} = 0.\]
Thus if $x,y \in I$, we can write $xy = (x_s+x_a)(y_s+y_a) = x_sy_s + x_sy_a + x_ay_s + x_ay_a$ and all four terms are taken to $0$, as desired: The first three terms all have a factor that is taken to zero by the calculation just made and $x_ay_a = (x_ay_a)_s$ and both $x_a$ and $y_a$ (and thus their product) are contained in $I$, so the same argument applies. It is hopefully clear that both maps are compatible with the respective structure maps from and to $S$.

Finally, we have to check that the two maps are inverse to one another: This is clear for the composition starting at 
$(S \otimes_R S)/I^2$, and for the composition starting at $(S \otimes_R S) \otimes_{\Delta_\varphi} S$ we compute 
\[\frac{st \otimes s' + s \otimes s't} 2 \otimes 1 = (s \otimes s' \otimes 1)\frac{1\otimes t \otimes 1+ t \otimes 1 \otimes 1} 2 \]\[= (s \otimes s' \otimes 1)(1 \otimes 1 \otimes t) = s \otimes s' \otimes t\]
by another application of (\ref{eq}).
\end{proof}

Identifying $I/I^2 = \Omega_\varphi$, the module of K\"ahler differentials, via the derivation $s \mapsto s \otimes 1 - 1 \otimes s$, and thus $\mathrm{Hom}_S(I/I^2,M) = \mathrm{Der}_R(S,M)$, we can unwind definitions in the splitting from Proposition \ref{prop:criterionintro} to find: 

\begin{corollary}
If $2$ is a unit in $S$, the map
\[\mathrm{Der}_R(S,M) \longrightarrow \Quad_\varphi(S^2,M), \quad d \longmapsto [(s,t) \mapsto d(s)t-sd(t)]\]
together with the inclusion $\Quad_S(S^2,M) \subseteq \Quad_\varphi(S^2,M)$ gives an isomorphism
\[\Quad_\varphi(S^2,M) = \Quad_S(S^2,M) \oplus \mathrm{Der}_R(S,M).\]
\end{corollary}

We obtain the first half of Theorem B:

\begin{corollary}
If $2$ is a unit in $S$, then the inclusion $\Quad_S(M,N) \subseteq \Quad_\varphi(M,N)$ is an equality for all $S$-modules $M$ and $N$ if and only if $I^2 = I$, where $I$ is the kernel of the multiplication $S \otimes_R S \rightarrow S$, or in other words if and only if $\Omega_\varphi = 0$.
\end{corollary}

At the other extreme we find:

\begin{proposition}\label{ident2=0}
  If $2=0$ in $S$, then we can regard $\mathrm{Fr}_\varphi^*S \otimes_{S \otimes_R \mathrm{Fr}^*R/2} \mathrm{Fr}_\varphi^*S$ as an $S$-algebra using \[S \longrightarrow \mathrm{Fr}_\varphi^*S \otimes_{S \otimes_R \mathrm{Fr}^*R/2} \mathrm{Fr}_{\varphi}^*S, \quad s \longmapsto (s \otimes 1),\]
where $\mathrm{Fr} \colon R \rightarrow R/2$ is the Frobenius of $R$ and $\mathrm{Fr}_{\varphi} \colon S \otimes_R \mathrm{Fr}^*R/2\rightarrow S$ is the relative Frobenius given by $s \otimes r \mapsto s^2\varphi(r)$.
Then the two maps
\[\QQ_\varphi \longrightarrow \mathrm{Fr}_\varphi^*S \otimes_{S \otimes_R \mathrm{Fr}^*R/2} \mathrm{Fr}_{\varphi}^*S, \quad s \otimes s' \otimes t \longmapsto t \otimes ss'\]
\[\mathrm{Fr}_\varphi^*S \otimes_{S \otimes_R \mathrm{Fr}^*R/2} \mathrm{Fr}_{\varphi}^*S  \longrightarrow \QQ_\varphi, \quad s \otimes s' \longmapsto s' \otimes 1 \otimes s\]
are inverse isomorphisms over and under $S$.
\end{proposition}

Note that $\mathrm{Fr}_\varphi^*S \otimes_{S \otimes_R \mathrm{Fr}^*R/2} \mathrm{Fr}_{\varphi}^*S$ is just a long name for $(S \otimes_R S)/J$, where $J$ is the ideal generated by the elements $st^2 \otimes s' - s \otimes t^2s'$ for $s,s',t \in S$. It can for example also be described as $S \otimes_{\Delta_\varphi} S$, but the description in terms of the relative Frobenius is what will allow us to deduce our main results.

\begin{proof}
The well-definedness of the first two maps is hopefully clear this time around. We showed in (\ref{eq}) that 
\[s \otimes s' \otimes 1 + s' \otimes s \otimes 1 = 1 \otimes 1 \otimes 2ss'\]
in $(S \otimes_R S) \otimes_{\Delta_\varphi} S$, which in the present case implies $s \otimes s' \otimes 1  = s' \otimes s \otimes 1$, and thus
\[s \otimes s' \otimes t = (s \otimes 1 \otimes 1)(1 \otimes s' \otimes 1)(1 \otimes 1 \otimes t) = (1 \otimes s \otimes 1))(1 \otimes s' \otimes 1)(1 \otimes 1 \otimes t) = 1 \otimes ss' \otimes t.\]
This implies the well-definedness of the last map via
\[s'u^2 \otimes 1 \otimes s = u \otimes s'u \otimes s = [(1 \otimes s')(u \otimes u)] \otimes s = 1 \otimes s' \otimes u^2s.\]
The compatibility with the structure maps to and from $S$ is hopefully evident.

The composition starting at $\mathrm{Fr}_\varphi^*S \otimes_{S \otimes_R \mathrm{Fr}^*R/2} \mathrm{Fr}_{\varphi}^*S$ is evidently the identity, and that starting at $(S \otimes_R S) \otimes_{\Delta_\varphi} S$ takes $s \otimes s' \otimes t$ to $ss' \otimes 1 \otimes t$ which is the identity by the (evident flip of the) relation above.
\end{proof}

Recalling that a ring map $R \rightarrow S$ is an epimorphism if and only if the multiplication map $S \otimes_R S \rightarrow S$ is an isomorphism, we deduce the second half of Theorem B:

\begin{corollary}
If $2 = 0$ in $S$, then the inclusion $\Quad_S(M,N) \subseteq \Quad_\varphi(M,N)$ is an equality for all $S$-modules $M$ and $N$ if and only if the map $\mathrm{Fr} \cdot \varphi \colon S \otimes R \rightarrow S$ is an epimorphism of rings.
\end{corollary}

We are not aware of a useful description of the complementary summand of $S$ inside $\mathrm{Fr}_\varphi^*S \otimes_{S \otimes_R \mathrm{Fr}^*R/2} \mathrm{Fr}_{\varphi}^*S$, and in particular, of its relation to Kähler differentials beyond the general surjection $\mathrm W_{S/R} \rightarrow \Omega_{S/R}$. \\

At any rate, with these characterisations established we can give the proof of our main result. 

\begin{proof}[Proof of Theorem A]
We start with the case of characteristic 2: Start by noting that an arbitrary ring map $\varphi \colon R \rightarrow S$ can be decomposed as $\varphi_3\varphi_2\varphi_1$, where $\varphi_1 \colon R \rightarrow R/\mathrm{ker}(\varphi)$ is a surjection, $\varphi_2 \colon R/\mathrm{ker}(\varphi) \rightarrow (R/\mathrm{ker}(\varphi))[U^{-1}]$ a localisation, where $U = \varphi^{-1}(S^\times)$, and finally $\varphi_3 \colon  (R/\mathrm{ker}(\varphi))[U^{-1}] \rightarrow S$ is injective and reflects units. Again, since both surjections and localisations are epimorphisms $\varphi$ is then an epimorphism if and only if $\varphi_3$ is one. But a map of this last type with target a field is then necessary a map between fields, and epimorphisms between fields are isomorphisms for dimension reasons (whenever $S$ is free over $R$ the multiplication map $S \otimes_R S \rightarrow S$ can clearly only be injective if the rank of $S$ is $1$). 

But this means that the map $\mathrm{Fr} \cdot \mathrm{incl} \colon L \otimes K \rightarrow L$ is an epimorphism, if and only if its image generates $L$ as a field, as desired. 

For the case of characteristic not $2$, we note that among fields the vanishing of the Kähler differentials is well-known to be equivalent to the conditions from the theorem; that the conditions given imply this vanishing is fairly simple: In characteristic $0$ it follows automatically that $L/K$ is separable algebraic and this generally implies the vanishing of $\Omega_{L/K}$: For any $a \in L$ the minimal polynomial $F \in K[T]$ of $a$ then has $F'(a) \neq 0$, so $0 = d(0) = d(F(a)) = F'(a)d(a)$, which implies $d(a) = 0$ and thus gives the conclusion, since the $d(a)$ form a generating set of $\Omega_{L/K}$ as $a$ ranges over $L$.

% as follows.
% Since everything is compatible with filtered colimits one reduces to the case where $L/K$ is finite. Either using the transitivity sequence for K\"ahler differentials for a chain of intermediate extension, or that every finite separable field extension admits a primitve element, we are thus reduced to considering the case $L \cong K[T]/f$, where $f$ is the minimal polynomial of some generator $x \in L$. We then find
%\[L \otimes_K L \cong L[T]/f\]
%with the action of $L \otimes_K L$-action on the right described by $x \otimes 1$ acting via multiplication with $x$, whereas $1 \otimes x$ acts by multiplication with $T$. But in $L[T]$ we have $f = (T-x)g$ with $g$ coprime to $T-x$ (by separability), so by the Chinese remainder theorem we find
%\[L \otimes_K L = L \times L[T]/g\]
%with the action of of $L \otimes_K L$ the obvious on the left factor, and the right factor still acted on as above. But in particular, the inclusion of $L$ into this defines an $L \otimes_K L$-linear section of the multiplication map $L \otimes_K L \rightarrow L$, and such ring maps always have idempotent kernels: If $\varphi \colon R \rightarrow S$ is a ring homomorphism which admits an $R$-linear section $s$, then one easily checks $(1-s(1))r = r$ for all $r \in \mathrm{ker}(\varphi)$, and clearly $1-s(1) \in \mathrm{ker}(\varphi)$, so the left hand side is in $\mathrm{ker}(\varphi)^2$.

In positive characteristic $p$ one immediately checks that all $K$-linear derivatives on $L$ are automatically linear over $L^p$, so in particular vanish if the extension $L/K$ is generated by $L^p$.

For the converses one needs a more detailed analysis: For characteristic $p>0$ (including $p=2$), a $p$-basis of the simple purely inseparable extension $L/K \cdot L^p$ gives a basis of $\Omega_{L/K \cdot L^p} \cong \Omega_{L/K}$ by applying the universal derivation, see e.g.\ \cite[Section 26.F]{Mat}. In characteristic $0$ a transcendence basis for $L/K$ does the same; this can be seen by combining the transitivity sequence \cite[Section 26.H]{Mat} for $L/E/K$, with $E$ generated by the transcendence basis, with the fact that $L \otimes_E \Omega_{E/K} \rightarrow \Omega_{L/K}$ is an isomorphism whenever $L/E$ is separable algebraic, which in turn follows from the conormal sequence \cite[Section 27.I]{Mat} as in \cite[Section 27.J]{Mat} reducing to the case of finite extensions by compatibility with filtered colimits, see for example also \cite[Korollar 7.4.9]{Bosch}.

In fact, the associated derivations detect both $p$- and transcendence bases according to the characteristic, but we do not need this.
\end{proof}

\section{Three examples}

Finally, we explicitly compute the universal rings $\QQ_{L/K}$ in three simple examples.  To this end, let $\Lambda_K$ denote exterior algebras over $K$ (with no signs inserted), i.e.\ 
\[\Lambda_K[S_1, \dots S_n] = K[S_1, \dots, S_n]/(S_i^2 \mid 1 \leq i \leq n).\]

\subsection*{Simple purely transcendental extensions}

Let us explicitly compute the universal ring for the field extension $K(T)/K$.

\begin{proposition}
The maps
\[\QQ_{K(T)/K} \longrightarrow K(T) \otimes_K \Lambda_K[S], \quad T \otimes 1 \otimes 1 \longmapsto T+S, \]\[ 1 \otimes T \otimes 1 \longmapsto T-S, \quad 1 \otimes 1 \otimes T \longmapsto T\]
\[K(T) \otimes_K \Lambda_K[S] \longrightarrow \QQ_{K(T)/K}, \quad T \longmapsto 1 \otimes 1 \otimes T, \]\[S \longmapsto T \otimes 1 \otimes 1 - 1 \otimes 1 \otimes T\]
are inverse isomorphisms, and so in particular $\mathrm{W}_{K(T)/K} \simeq K(T)$.

Consequently, the $K(T)$-vector space 
\[\mathrm{ker}\left(\Quad_{K(T)/K}(K(T)^2,K(T)) \rightarrow \Quad_{K(T)}(K(T),K(T))^2\right)\]
of $K(T)/K$-quadratic forms on $K(T)^2$, that vanish on both $K(T) \oplus 0$ and $0 \oplus K(T)$, is isomorphic to 
\[\mathrm{Hom}_{K(T)}(K(T) \otimes_K \Lambda_K[S],K(T)) = \mathrm{Hom}_K(\Lambda[S],K(T)) = K(T)^2,\] with $f \colon K(T) \otimes_K \Lambda_K[S] \rightarrow K(T)$ corresponding to the form
\[(F,G) \longmapsto f(F(T+S)G(T-S)) = FGf(1)+(F'G-FG')f(S)\]
and similarly
\[\Quad_{K(T)/K}(K(T)^2,K(T)) = \Quad_{K(T)}(K(T)^2,K(T)) \oplus K(T)\]
\end{proposition}

If the characteristic of $K$ is not $2$ then this follows immediately from the results of the previous section by observing $\Omega_{K(T)/K} = K(T)$, but the argument we give below works simultaneously in characteristic $2$.

\begin{proof}
Immediately from the definition and Lemma \ref{lemma:flatfixed} we find that $\QQ_{K(T)/K}$ is the localisation of \[K[X,Y] \otimes_{K[X+Y,XY]} K[T]\]
at all non-zero polynomials in $X, Y$ or $T$; here, the tensor product is formed along the homomorphism determined by $X+Y \mapsto 2T$ and $XY \mapsto T^2$. \\

This ring (before localisation) is easily checked isomorphic to $K[T] \otimes_K \Lambda_K[S]$ via 
\[K[T] \otimes_K \Lambda_K[S] \longrightarrow K[X,Y] \otimes_{K[X+Y,XY]} K[T], \quad T \longmapsto T, S \longmapsto X-T\]
\[K[X,Y] \otimes_{K[X+Y,XY]} K[T] \longrightarrow K[T] \otimes_K \Lambda_K[S], \quad X \longmapsto T+S, Y \longmapsto T-S, T \longmapsto T.\]
This translates the structure maps to the maps induced by
\[K[T] \longrightarrow K[T] \otimes_K \Lambda_K[S]  \quad T \longrightarrow T\]
\[K[T] \otimes_K \Lambda_K[S] \longrightarrow K[T], \quad T\longrightarrow T, S \longrightarrow 0\]
and the universal form is induced by
\[K[T] \otimes_K K[T] \longrightarrow K[T] \otimes_K \Lambda_K[S], F \otimes G \longmapsto F(T+S)G(T-S).\]
  The universal ring is then the localisation at all non-zero polynomials in $T,T+S,T-S$, but after inverting the former to get $K(T) \otimes_K \Lambda_K[S]$ all non-zero polynomials in $T+S$ and $T-S$ have already become invertible, since 
\[F(T+S) = F(T)+F'(T)S = F(T)\left(1+\frac{F'(T)}{F(T)} S\right)\]
\[F(T-S) = F(T)-F'(T)S = F(T)\left(1-\frac{F'(T)}{F(T)} S\right)\]
the product of two units in both cases. This shows that indeed $\QQ_{K(T)/K} \cong K(T) \otimes_K \Lambda_K[S]$ and that a homomorphism $f \colon K(T) \otimes_K \Lambda[S] \rightarrow K(T)$ corresponds to the form sending $(F,G) \in K(T)^2$ to 
\begin{align*}
f(F(T+S)G(T-S)) &= f(F(T)+F'(T)S)(G(T)-G'(T)S)) \\
&= FGf(1)+(F'G-FG')f(S)
\end{align*}
as claimed.

To finally obtain the claimed description of 
\[\mathrm{ker}\left(\Quad_{K(T)/K}(K(T)^2,K(T)) \rightarrow \Quad_{K(T)}(K(T),K(T))^2\right)\]
simply its identification with $\mathrm{Hom}_{K(T)}(\QQ_{K(T)/K}, K(T))$ from Lemma \ref{lm:quadsplit}.
\end{proof}

\subsection*{Simple purely inseparable extensions}

Let us explicitly compute the universal ring for the field extension $K(T)/K(T^p)$ when $\mathrm{char}(K) = p$; the argument is very similar to the previous case, and in particular, the result is the same:

\begin{proposition}
If $K$ has characteristic $p$, the maps
\[\QQ_{K(T)/K(T^p)} \longrightarrow K(T) \otimes_K \Lambda_K[S], \quad T \otimes 1 \otimes 1 \longmapsto T+S, \]\[ 1 \otimes T \otimes 1 \longmapsto T-S, \quad 1 \otimes 1 \otimes T \longmapsto T\]
\[K(T) \otimes_K \Lambda_K[S] \longrightarrow \QQ_{K(T)/K(T^p)}, \quad T \longmapsto 1 \otimes 1 \otimes T, \]\[S \longmapsto T \otimes 1 \otimes 1 - 1 \otimes 1 \otimes T\]
are inverse isomorphisms, so in particular $\mathrm{W}_{K(T)/K(T^p)} \cong K(T)$.

Consequently, the $K(T)$-vector space 
\[\mathrm{ker}\left(\Quad_{K(T)/K(T^p)}(K(T)^2,K(T)) \rightarrow \Quad_{K(T)}(K(T),K(T))^2\right)\]
of $K(T)/K(T^p)$-quadratic forms on $K(T)^2$, that vanish on both $K(T) \oplus 0$ and $0 \oplus K(T)$, is isomorphic to 
\[\mathrm{Hom}_{K(T)}(K(T) \otimes_K \Lambda[S],K(T)) = \mathrm{Hom}_K(\Lambda[S],K(T)) = K(T)^2,\] with $f \colon K(T) \otimes_K \Lambda[S] \rightarrow K(T)$ corresponding to the form
\[(F,G) \longmapsto f(F(T+S)G(T-S)) = FGf(1)+(F'G-FG')f(S).\]
and similarly
\[\Quad_{K(T)/K(T^p)}(K(T)^2,K(T)) = \Quad_{K(T)}(K(T)^2,K(T)) \oplus K(T)\]
\end{proposition} 

\begin{proof}
Immediately from the definition and Lemma \ref{lemma:flatfixed} we find that $\QQ_{K(T)/K(T^p)}$ is the localisation of \[K[X,Y]/(X^p-Y^p) \otimes_{K[X+Y,XY]} K[T]\]
at all non-zero polynomials in $X, Y$ or $T$. But then performing the same steps as in the previous proof, we find that this is $(K[T] \otimes_K \Lambda_K[S])/((T-S)^p - (T+S)^p)$ localised at all non-zero polynomials in $T$. But $(T-S)^p - (T+S)^p = -2S^p = 0$ already in $K[T] \otimes_K \Lambda_K[S]$, so this does not affect the proof any further.
\end{proof}

In both cases, we see that even in characteristic $2$ the examples of non-trivial $L/K$-quadratic forms come from $K$-linear derivations on $L$.

\subsection*{General purely transcendental extensions}
The uniform behaviour exhibited by the examples above disappears in the case of more than one variable. Proposition \ref{ident2=0} implies: 

\begin{proposition}
If the characteristic of $K$ is $2$, then the maps
\[\QQ_{K(T_1, \dots T_n)/K} \longrightarrow K(T_1, \dots, T_n) \otimes \Lambda_K[S_1, \dots, S_n], \]
\[T_i \otimes 1 \otimes 1 \longmapsto T_i \otimes 1 + 1 \otimes S_i, \quad 1 \otimes T_i \otimes 1 \longmapsto  T_i \otimes 1 + 1 \otimes S_i,\]
\[1 \otimes 1 \otimes T_i \longmapsto T_i \otimes 1\]
and
\[K(T_1, \dots, T_n) \otimes \Lambda_K[S_1, \dots, S_n] \longrightarrow\QQ_{K(T_1, \dots T_n)/K}\]
\[T_i \otimes 1 \longmapsto 1 \otimes 1 \otimes T_i, \quad 1 \otimes S_i \longmapsto 1 \otimes 1 \otimes T_i + T_i \otimes 1 \otimes 1\]
are inverse isomorphisms. In particular, the dimension of $\mathrm{W}_{K(T_1, \dots T_n)/K}$ is $2^n-1$ and the map $W_{K(T_1, \dots T_n)/K} \longrightarrow \Omega_{K(T_1, \dots T_n)/K}$, whose target has dimension $n$, is not injective for $n\geq 2$; it vanishes on the square of the augmentation ideal $I$ of $\QQ_{K(T_1, \dots T_n)/K}$.
\end{proposition}

Explicitly, denoting by $\partial^k(-)/\partial T_{i_1}\cdots\partial T_{i_k} \colon K(T_1, \dots, T_n) \rightarrow K(T_1, \dots T_n)$ the $k$-fold derivative in the directions of $T_{i_1}, \dots, T_{i_k}$, we find the additional $K(T_1, \dots T_n)/K$-quadratic forms $K(T_1, \dots T_n)^2 \rightarrow K(T_1, \dots T_n)$ spanned by the maps
\[q_{i_1, \dots i_k}(F,G) = \frac{\partial^k FG}{\partial T_{i_1}\cdots\partial T_{i_k}}.\]
for monotone injections $i \colon \{1,\dots,k\} \rightarrow \{1, \dots, n\}$. Note that the unique such $i$ for $k=0$ yields the sole $K(T_1, \dots T_n)$-quadratic form, that vanishes on $K(T_1, \dots T_n) \times 0$ and $0 \times K(T_1, \dots T_n)$, whereas the $n$ maps with $k=1$ span the $n$-dimensional subspace of forms coming from derivations as in the introduction; in characteristic $2$ we have $d(F)G-Fd(G) = d(FG)$ for any derivation $d$ after all.

The forms parametrised by injections with $k\geq 2$ are, however, new and specific to characteristic $2$: If the characteristic of $K$ is not $2$, then the map 
\[\QQ_{K(T_1, \dots T_n)/K} \longrightarrow K(T_1, \dots, T_n) \otimes \Lambda_K[S_1, \dots, S_n]/I^2, \]
is an isomorphism by \ref{ident2inv}, so that the forms involving multiple derivatives have no analogue away from characteristic $2$.

\begin{proof}
Using \ref{ident2=0} we find
\[\QQ_{K[T_1, \dots T_n]/K} \cong K[T_1, \dots, T_n] \otimes_{K[T_1^2, \dots, T_n^2]} K[T_1, \dots, T_n]\]
In the latter ring we have 
\[(1 \otimes T_i + T_i \otimes 1)^2 = 1 \otimes T_i^2 + T_i^2 \otimes 1 = 2 \otimes T_i^2 = 0\] 
which makes the map
\[K[T_1, \dots, T_n] \otimes_K \Lambda_K[S_1, \dots, S_n] \longrightarrow K[T_1, \dots, T_n] \otimes_{K[T_1^2, \dots, T_n^2]} K[T_1, \dots, T_n]\]
\[T_i \longmapsto T_i \otimes 1, S_i \longmapsto 1 \otimes T_i + T_i \otimes 1\] 
well-defined. Since it is obviously surjective and both sides are free of dimension $2^n$ over $K[T_1, \dots, T_n]$ it is an isomorphism. Localising this isomorphism just as in the previous proofs gives the result.

For the final claim, recall that the map $W_{S/R} \rightarrow \Omega_{S/R}$ is given by $s \otimes s' \otimes t \mapsto t(d(s)s'-sd(s'))$. Under the identifications above this map then takes $1 \otimes S_i$ to $d(T_i)$. For products of these elements we first note that the identification above takes $\prod_{j=1}^k 1 \otimes S_{i_j}$ to $\sum_{P \subseteq \{1,\dots, k\}} (\prod_{p \in P} T_{i_p}) \otimes (\prod_{p \notin P} T_{i_p})$. This is in turn taken to
\[\sum_{P \subseteq \{1,\dots, k\}} \left(\prod_{p \in P} T_{i_p}\right) \cdot d\left(\prod_{p \notin P} T_{i_p}\right) = \sum_{{P \subseteq \{1,\dots, k\}} \atop{j \notin P}} \left(\prod_{p \neq j} T_{i_p}\right) d(T_{i_j})\]\[= 2^{k-1} \sum_{j = 1}^k  \left(\prod_{p \neq j} T_{i_p}\right)d(T_{i_j})\]
which vanishes for $k >1$. 
\end{proof}

\bibliographystyle{amsalpha}
\newcommand{\etalchar}[1]{$^{+}$}

\end{document}